\documentclass{article}

\usepackage{amsmath}
\usepackage{amsfonts}
\usepackage{amssymb}
\usepackage{amsthm}
\usepackage{array}
\usepackage{chngpage}
\usepackage{comment}
\usepackage{float}
\usepackage[OT2,T1]{fontenc}
\usepackage{graphicx,caption}
\usepackage[export]{adjustbox}
\usepackage{longtable}
\usepackage{mathrsfs}
\usepackage{mathtools}
\usepackage{wrapfig}
\usepackage{cutwin}
\usepackage{shapepar}
\usepackage{tikz}
\usepackage[lowtilde]{url}
\usepackage{subcaption}
\usetikzlibrary{matrix,arrows}
\usepackage{tabularx}
\usepackage{multirow}

\newcommand{\C}{\mathbb{C}}

\newcommand{\F}{\mathbb{F}}
\newcommand{\HH}{\mathbb{H}}

\newcommand{\Q}{\mathbb{Q}}
\newcommand{\R}{\mathbb{R}}
\newcommand{\Z}{\mathbb{Z}}

\newcommand{\cC}{\mathcal{C}}

\newcommand{\cH}{\mathcal{H}}
\newcommand{\cM}{\mathcal{M}}
\newcommand{\cO}{\mathcal{O}}

\DeclareSymbolFont{cyrletters}{OT2}{wncyr}{m}{n}
\DeclareMathSymbol{\sha}{\mathalpha}{cyrletters}{"58}

\newcommand{\eps}{\varepsilon}

\newlength{\strutheight}
\newcommand{\thalf}{\tfrac{1}{2}}

\newtheorem{theorem}{Theorem}[section]

\newtheorem{proposition}[theorem]{Proposition}

\author{Alex Cowan\\cowan@math.harvard.edu}
\title{Computing newforms using supersingular isogeny graphs}
\date{}

\begin{document}
\maketitle
\begin{abstract}
  \noindent
  We describe an algorithm that we used to compute the $q$-expansions of all weight $2$ cusp forms of prime level at most $2,\!000,\!000$ and dimension at most $6$. We also present an algorithm that we used to verify that there was only one cusp form of dimension $7$ or more per Atkin-Lehner eigenspace for prime levels between $10,\!000$ and $1,\!000,\!000$. Our algorithm is based on Mestre's M\'{e}thode des Graphes, and involves supersingular isogeny graphs and Wiedemann's algorithm for finding the minimal polynomial of sparse matrices over finite fields.
\end{abstract}
\tableofcontents

\section{Introduction}\label{intro_section}
Let $S_2(p)$ be the space of weight $2$ cusp forms of prime level $p$. We say that the \textit{dimension} of a newform of $S_2(p)$ is the degree of the number field its Hecke eigenvalues generate, or, equivalently, the size of its Galois orbit. There have been many efforts from computational number theorists to create databases with information about the newforms of $S_2(p)$; we highlight the Antwerp tables \cite{antwerp}, Cremona's database of elliptic curves \cite{cremona}, and the LMFDB \cite{lmfdb}. The LMFDB builds on the previous two, and currently lists the $q$-expansions of every newform of $S_2(p)$ of dimension $g$ at most $20$ and level at most $10,\! 000$ \cite{bbbccdldlrsv}.\\
\\
In this paper, we describe an algorithm that we used to compute the $q$-expansions of all newforms with $g \leq 6$ and $p < 2,\!000,\!000$ up to the Sturm Bound, and also used to verify that there were no eigenforms with $g \geq 7$ and $10,\!000 < p < 1,\!000,\!000$ besides one factor of high dimension per Atkin-Lehner eigenspace. For each level $p$, our algorithm runs in time $\cO(p^{2+\eps})$ and space $\cO(p^{1+\eps})$.\\
\\
In the field of arithmetic statistics there is a lot of interest in understanding how various properties of these eigenforms are distributed. This is partly because modular forms are interesting in their own right, but also partly because standard modularity conjectures \cite{diamond-shurman} predict that, for each genus $g$ factor of the modular jacobian $J_0(p)$, there is an associated weight $2$ newform of level $p$ and dimension $g$. The genus $1$ case of elliptic curves has been studied extensively, and is one of the most important topics in modern number theory. The association between elliptic curves and $1$-dimensional modular forms is a theorem \cite{wiles, taylor-wiles, bcdt}, and the literature contains conjectures and theorems for how many related invariants are distributed, notably ranks \cite{ppvw, bklhpr}, Selmer groups \cite{bhargava-shankar}, torsion subgroups \cite{harron-snowden}, and other numbers which appear in the Birch and Swinnerton-Dyer formula \cite{poonen}. Generalizations of these theorems to genus $2$ or more are far out of reach for the most part, and in many situations it is poorly understood what the correct generalizations would even be. In particular, merely predicting the number of genus $g$ factors of $J_0(p)$ has not been done whenever $g \geq 2$, whereas there are well established conjectures for the number of elliptic curves with bounded conductor \cite{brumer-mcguinness, watkins}. In light of this gap in understanding, databases of newforms of $S_2(p)$ are very useful: they give a way to observe generalizations of phenomena which occur in the genus $1$ case, and they also allow one to formulate conjectures about the statistics of these objects. In \cite{martin}, for instance, Martin computes the dimensions of the eigenforms of $S_2(p)$ for $p < 60000$ and uses this data to formulate conjectures related to counts of eigenforms of fixed dimension.\\
\\
The main idea in our algorithm is from Mestre's M\'{e}thode des Graphes \cite{mestre}. In Mestre's work, he relates the $q$-expansion of weight $2$ newforms of prime level to ``supersingular isogeny graphs''. The \textit{supersingular $\ell$-isogeny graph $\text{over }\overline{\F}_p$} is the graph whose vertices are supersingular $j$-invariants over $\overline{\F}_p$, and whose edges are $\ell$-isogenies. These graphs have recently been of independent interest because of their applications to cryptography \cite{charles-lauter-goren, jao-defeo, ehlmp}. The relationship Mestre highlights depends on a trace formula: the action of the Hecke operator $T_\ell$ on the space $S_2(p)$ can be represented as the adjacency matrix of the supersingular $\ell$-isogeny graph. We find simultaneous eigenvectors of these matrices, and then use a formula from Mestre's work to compute the associated $q$-expansions.\\
\\
Some of the building blocks of our algorithm come from more general-purpose techniques. In particular, we use Wiedemann's algorithm for finding the minimal polynomial of a sparse matrix over a finite field \cite{wiedemann}, Brent and Kung's algorithm for fast power series composition \cite{brent-kung}, and a well known dynamic programming algorithm for solving the subset sum problem \cite{clrs}.\\
\\
We implemented most of the algorithm in Sage. We used cython to multiply vectors by sparse matrices, and we used \texttt{c++}'s NTL package when doing power series manipulations. We ran our code on the Oklahoma University Supercomputing Center for Education \& Research.

\section{Acknowledgments}\label{acknowledgements_section}
I'm extremely grateful to Noam Elkies for both suggesting I work on this problem, and for telling me about many of the ideas and techniques used in this paper. I thank Kimball Martin for many helpful discussions regarding ways to use the data, and for running and managing the computations. The computing for this project was performed at the OU Supercomputing Center for Education \& Research (OSCER) at the University of Oklahoma (OU). I thank Drew Sutherland for help in adding the data to the LMFDB. I thank David Roe for helping me write some of the code. This work was done under Simons Collaboration grant number 550031.

\section{Background}\label{background_section}
\subsection{Wiedemann's algorithm}\label{background_wiedemann_subsection}
We use \cite{wiedemann} as a reference for what we call ``Wiedemann's algorithm''. Given an $n\times n$ nonsingular matrix $M$ over a finite field $F$, Wiedemann gives a probabilistic algorithm for finding the minimal polynomial $\mu$ of $M$. Wiedemann's algorithm is one of the major building blocks of the algorithm presented in this paper. See section \ref{charpoly_section} for details.\\
\\
Let $u$ be a vector in $F^n$, and let $i \in [1,n]\cap \Z$ be one of the indices of the coordinates of $u$. The key idea in Wiedemann's algorithm is that the sequence $u_i, (Mu)_i, (M^2u)_i, \dots$ will satisfy a recursion relation, and that recursion relation will, for most choices of $u$ and $i$, give the minimal polynomial $\mu$.\\
\\
For a given choice of $u$ and $i$, computing the sequence $u_i, (Mu)_i, (M^2u)_i, \dots, (M^ru)_i$ takes time $\cO(r\omega)$, where $\omega$ is the number of nonzero entries in $M$. Wiedemann then uses the following:
\begin{proposition}\label{wiedemann_recurrence}
  If $\mu(t) = \sum_{k=0}^n \mu_n t^n$ is the minimal polynomial of $M$, then, for any vector $u$, we have
  $$\sum_{k=0}^n \mu_n M^k u = 0.$$
\end{proposition}
\noindent
As a consequence of proposition \ref{wiedemann_recurrence}, we know that the sequence $u_i, (Mu)_i, (M^2u)_i, \dots, (M^ru)_i$ satisfies a recursion relation of length at most $n$. To determine what the recursion relation is, we'll need a number of terms at least double the recursion length. Thus, we can take $r = 2n$. In our application, $M$ will have $\cO(n)$ nonzero entries, so we'll compute the sequence $u_i, (Mu)_i, (M^2u)_i, \dots, (M^{2n}u)_i$ in time $\cO(n^2)$.\\
\\
Wiedemann then uses an algorithm of Berlekamp-Massey to find the coefficients of the recursion relation of the sequence $u_i, (Mu)_i, (M^2u)_i, \dots, (M^{2n}u)_i$. This takes time $\cO(n^2)$. See \cite{mills} for a description of the Berlekamp-Massey algorithm in terms of continued fractions. One step in this algorithm involves writing $\frac{1}{\mu(t)}$ as a power series. If $M$ is not invertible, then $\mu(t)$ will be divisible by $t$, making it impossible to do this. It's possible to circumvent this problem in a number of ways. We chose to modify our matrices so that they would be invertible (see section \ref{charpoly_shift_subsection}), but Wiedemann gives a modification of his algorithm for this case in \cite{wiedemann}, and it's also possible to modify the Berlekamp-Massey algorithm directly.
\subsection{Newforms}\label{background_modforms_subsection}
We give \cite{darmon} as a reference for this section. Define $\Gamma_0(p)$ to be the group
$$\Gamma_0(p) \coloneqq \left\{\begin{pmatrix}a & b\\c & d\end{pmatrix} \in \text{SL}_2(\Z) \,:\, c = 0\text{ mod } p\right\}.$$
The group $\Gamma_0(p)$ acts on the upper half-plane $\HH \coloneqq \{x + iy \,:\, x \in \R, y \in \R_{>0}\}$ via M\"{o}bius transformations:
$$\begin{pmatrix}a & b\\c & d\end{pmatrix} z \coloneqq \frac{az + b}{cz + d}.$$
A \textit{weight $2$ modular form on $\Gamma_0(p)$} is a holomorphic function $f : \HH \to \C$ which satisfies the relation
$$f\!\left(\begin{pmatrix}a & b\\c & d\end{pmatrix} z\right) = (cz + d)^2 f(z)$$
for every $z \in \HH$ and every $\begin{pmatrix}a & b\\c & d\end{pmatrix} \in \Gamma_0(p)$. Weight $2$ modular forms on $\Gamma_0(p)$ form a finite-dimensional complex vector space.\\
\\
Because $\begin{pmatrix} 1 & 1\\0 & 1\end{pmatrix}$ is in $\Gamma_0(p)$, we have $f(z + 1) = f(z)$. It follows that modular forms have Fourier expansions, i.e. there exist complex numbers $a_n(f)$ such that
$$f(z) = \sum_{n = 0}^\infty a_n(f) e^{2\pi i n z}.$$
The space $S_2(p)$ of cusp forms is the subspace of these modular forms which have $a_0 = 0$. Throughout this paper we use the shorthand $q \coloneqq e^{2\pi i z}$, and we'll call these Fourier expansions \textit{$q$-expansions}. We'll write $a_n$ instead of $a_n(f)$ when the modular form $f$ is clear from context.\\
\\
The \textit{Hecke operators} are linear operators $T_n$ indexed by positive integers $n$ which act on the space of modular forms as
$$a_m(T_nf) = \sum_{d | \text{gcd}(m,n)} d \cdot a_{mn/d^2}(f).$$
The Hecke operators commute with one another, and hence they are simultaneously diagonalizable. Thus, there exist modular forms $f$ with $a_1 = 1$ which, for all $n$ simultaneously, satisfy the relations
$$T_n f = a_n f.$$
Modular forms with these properties are called \textit{newforms}. They form a basis for the space $S_2(p)$.\\
\\
In \cite{sturm}, Sturm proves the following:
\begin{theorem}[Sturm bound] If $f$ and $g$ are weight $2$ newforms of level $p$ and $a_n(f) = a_n(g)$ for all $n \leq \lfloor \frac{p+1}{6}\rfloor$, then $f = g$.
\end{theorem}
\noindent
There's are analogous results for other weights and levels as well.\\
\\
There is a linear operator $w_p$ called the \textit{Atkin-Lehner involution} which acts on $S_2(p)$. This operator commutes with all of the Hecke operators, so newforms are also eigenforms of the Atkin-Lehner involution. As suggested by the name, the Atkin-Lehner involution is an involution. Thus, if $f$ is a newform, then $w_p f = \pm f$, and $S_2(p)$ decomposes into two Atkin-Lehner eigenspaces (of roughly equal size; see \cite{martin_alsize}).

\subsection{La M\'{e}thode des Graphes}\label{background_ssj_subsection}
The \textit{supersingular $j$-invariants over $\overline{\F}_p$} are the $j$-invariants of the elliptic curves defined over $\overline{\F}_p$ which are supersingular (i.e. their endomorphism ring is an order in a quaternion algebra). It's known that there are $\lfloor\frac{p}{12}\rfloor + 0, 1, \text{or }2$ supersingular $j$-invariants over $\overline{\F}_p$, and that they're all defined over $\F_{p^2}$ \cite{silverman}. The \textit{supersingular $\ell$-isogeny graph over $\overline{\F}_p$} is the directed multigraph whose vertices are the supersingular $j$-invariants over $\overline{\F}_p$, and whose edges correspond to $\ell$-isogenies over $\overline{\F}_p$ between the associated elliptic curves. These graphs have been and continue to be studied extensively \cite{kohel, sutherland, charles-lauter-goren, jao-defeo, ehlmp}, in part because of their potential applications to post-quantum cryptography.\\
\\
As described in \cite[$\S$2,5]{gross}, \cite{mestre}, and \cite{dechene}, one description of the action of the Hecke operator $T_\ell$ on $S_2(p)$ is as the adjacency matrix of the supersingular $\ell$-isogeny graph over $\overline{\F}_p$. This connection involves a trace formula and an equivalence of categories between supersingular elliptic curves and orders in quaternion algebras.\\
\\
As a consequence of the representation of $T_\ell$ as the adjacency matrix of the supersingular $\ell$-isogeny graph over $\overline{\F}_p$, there is a bijection between newforms $f$ of $S_2(p)$ and vectors $v = (v_j)_j$ with coordinates indexed by the supersingular $j$-invariants over $\overline{\F}_p$ which are simultaneous eigenvectors of all the Hecke operators $T_\ell$. Mestre \cite{mestre} uses this bijection to produce the identity of power series
\begin{equation*}\label{mestre_eqn_background}
  \left(\sum_j v_jj\right)f(q)\frac{dq}{q} = \sum_j v_j\frac{dj(q)}{j(q) - j} \quad\text{ mod }\mathfrak{p},
\end{equation*}
where
\begin{itemize}
\item $j(q)$ is the modular $j$ function,
\item $\mathfrak{p}$ is any prime above $p$ in the number field generated by the Hecke eigenvalues $a_\ell$, and
\item the sums are over the supersingular $j$-invariants over $\overline{\F}_p$.
\end{itemize}
\noindent
The Weil bound for $f$ states that $|a_n| < 2\sqrt{n}$, so for $n \ll p^2$ this equality of power series is enough to know the values of $a_n$ exactly. Since the Sturm bound is $\cO(p)$, this identity is enough to distinguish newforms.

\section{Overview of the algorithm}\label{alg_section}
\noindent
The algorithm is broken into five sections:
\begin{itemize}
\item Section \ref{T2_section}: Computing the action of $T_2$ on $S_2(p)$.
\item Section \ref{charpoly_section}: Computing the characteristic polynomial of $T_2$ modulo some small auxiliary prime.
\item Section \ref{lift_section}: Determining the eigenvalues of degree $6$ or less $T_2$, and finding eigenbases over $\Z$ whenever they'll correspond to newforms of dimension $6$ or less.
\item Section \ref{qexpansion_section}: Computing $q$-expansions using Mestre's formula.
\item Section \ref{highdeg_section}: Verifying that only one high genus factor per Atkin-Lehner eigenspace exists.
\end{itemize}
\noindent
In section \ref{T2_section}, we compute a representation of the action of $T_2$ on $S_2(p)$, or, more precisely, its action on each Atkin-Lehner eigenspace, by constructing a supersingular isogeny graph. We use ``modular polynomials'' to find edges, and explore the graph using a breath first search.\\
\\
In section \ref{charpoly_section}, we compute the characteristic polynomial $\chi_\nu$ of $T_\ell$ modulo some small auxiliary prime $\nu$. The main ingredient in this step is Wiedemann's algorithm \cite{wiedemann} for computing minimal polynomials of sparse matrices defined over finite fields. This part of the algorithm contributes to the leading term in the overall asymptotic time complexity.\\
\\
In section\ref{lift_section}, we use the characteristic polynomial $\chi_\nu$ to find the eigenvalues and the eigenvectors of $T_\ell$ which correspond to low degree factors of its characteristic polynomial $\chi_\Z$ over $\Z$. Our method for ``lifting eigenspaces'' in this way is based on the heuristic that each low degree factor of $\chi_\Z$ has an associated eigenbasis made up of vectors whose coordinates are small.\\
\\
In section \ref{qexpansion_section}, we use a formula from Mestre's M\'{e}thode des Graphes to get the $q$-expansion of the newforms in terms of a power series involving the previously computed eigenvectors. In evaluating this formula, we used an algorithm from \cite{brent-kung} for composing power series.\\
\\
Finally, in section \ref{highdeg_section}, we find the degrees of the irreducible factors of $\chi_\Z$. This allows us to know how $J_0(p)$ decomposes as a product of Abelian varieties. As we mentioned before, we verified that, for $10^4 < p < 10^6$, there was only one factor per Atkin-Lehner eigenspace which had dimension $7$ or more. This part of the algorithm uses a modified version of a well known dynamic programming algorithm for solving the subset-sum problem \cite{clrs}. It is also a leading term in the asymptotic time complexity, and in practice is the most computationally expensive part of our algorithm, but it is also optional: the $q$-expansions of all newforms associated to factors of $J_0(p)$ of dimension at most $6$ can be computed independently of this verification that only two other irreducible factors per level exist.

\section{Computing the action of $T_\ell$ on $S_2(p)$}\label{T2_section}
For a given prime $\ell$, we generate two directed weighted multigraphs, $G_\ell^+$ and $G_\ell^-$, whose adjacency matrices are representations of the action of $T_\ell$ on the $+$ and $-$ Atkin-Lehner eigenspaces of $S_2(p)$ respectively.\\
\\
The supersingular $j$-invariants over $\overline{\F}_p$ are all defined over $\F_{p^2}$, and it's convenient to pick a generator of $\F_{p^2}$ which has trace $0$ because it simplifies the part of the algorithm described in section \ref{qexpansion_section}. Let $\sigma$ be the nontrivial element of $\text{Gal}(\F_{p^2}/\F_p)$. Pick some arbitrary ordering $<$ of the supersingular $j$-invariants over $\overline{\F}_p$. The vertices of $G_\ell^+$ are the pairs $(j,j^\sigma)$ with $j < j^\sigma$, and the vertices of $G_\ell^-$ are the pairs $(j,j^\sigma)$ with $j \leq j^\sigma$.\\
\\
The graph $G_\ell^+$ has a weight $1$ edge from $(j_1, j_1^\sigma)$ to $(j_2, j_2^\sigma)$ for each $\ell$-isogeny from $j_1$ to $j_2$ and for each $\ell$-isogeny from $j_1^\sigma$ to $j_2^\sigma$, and has a weight $-1$ edge from $(j_1, j_1^\sigma)$ to $(j_2, j_2^\sigma)$ for each $\ell$-isogeny from $j_1$ to $j_2^\sigma$ and for each $\ell$-isogeny from $j_1^\sigma$ to $j_2$.\\
\\
The graph $G_\ell^-$ has a weight $1$ edge from $(j_1, j_1^\sigma)$ to $(j_2, j_2^\sigma)$ for each $\ell$-isogeny from $j_1$ to $j_2$, $j_1$ to $j_2^\sigma$, $j_1^\sigma$ to $j_2$, and $j_1^\sigma$ to $j_2^\sigma$.\\
\\
Constructing these graphs is done in two steps: finding a starting vertex (section \ref{T2_startingvertex_subsection}), and exploring the graph (section \ref{T2_exploring_subsection}).

\subsection{Finding a starting vertex}\label{T2_startingvertex_subsection}
Let $j \in \Z$ be the $j$-invariant of an elliptic curve over $\Q$ with complex multiplication, and let $D$ be the discriminant of the associated imaginary quadratic field. The reduction $j$ mod $p$ is a supersingular $j$-invariant if $D$ is not a square mod $p$ \cite{broker}. Thus, over $99\%$ of the time, the reduction of one of the $13$ supersingular $j$-invariants over $\Q$ will be a supersingular $j$-invariant over $\overline{\F}_p$. We use this as our starting vertex in these cases.\\
\\
If $p$ is such that every $D$ is a square mod $p$, then we use code from Arpin, Camacho-Navarro, Lauter, Lim, Nelson, Scholl, and Sot{\'a}kov{\'a} to get the starting vertex \cite{adventures}.

\subsection{Exploring the graph}\label{T2_exploring_subsection}
For each prime $\ell$ there is a \textit{modular polynomial} $\phi_\ell(x,y) \in \Z[x,y]$ with the property that $\phi_\ell(j,y)$ has a zero at $y = j'$ of order equal to the number of $\ell$-isogenies from $j$ to $j'$ \cite{sutherland}. To generate the graphs $G_\ell^+$ and $G_\ell^-$, we do a breadth first search, finding the roots of $\phi_\ell(j,y)$ at the vertex $(j,j^\sigma)$ at each step. Because $\phi_\ell$ has coefficients in $\Z$, the roots of $\phi_\ell(j^\sigma, y)$ are the Galois conjugates of the roots of $\phi_\ell(j,y)$. We also make use of the fact that, beyond the first vertex, we know at least one of the roots of $\phi_\ell(j,y)$, reducing the degree of the polynomial we have to solve by $1$. Thus generating $G_\ell^\pm$ requires finding the roots in $\F_{p^2}$ of one polynomial of degree $\ell+1$ and $\cO(p)$ polynomials of degree $\ell.$\\
\\
In our implementation we take $\ell = 2$, but later in the algorithm it is sometimes necessary to compute the action of $T_\ell$ for $\ell \geq 3$. We do this without making use of any information gained while computing the action of $T_2$, because in practice this ran the fastest. We tested a different algorithm which made use of the fact that if there's an $\ell_1$-isogeny from $j$ to $j_1$, and an $\ell_2$-isogeny from $j$ to $j_2$, then there necessarily exists a $j'$ which is both $\ell_1$-isogenous to $j_2$ and $\ell_2$-isogenous to $j_1$, but our implementation took longer to evaluate the corresponding modular polynomials at all the neighbours of $j_1$ and $j_2$ than to find the roots of the modular polynomials directly.

\section{Computing the characteristic polynomial of $T_2$ mod $\nu$}\label{charpoly_section}
The next step of our algorithm for computing $q$-expansions is computing the characteristic polynomials of the adjacency matrices of the graphs $G_\ell^\pm$ (which we'll denote $T_\ell$ in a slight abuse of notation). Computing the characteristic polynomial over $\Z$ directly appears to be infeasible. Instead we compute modulo some small arbitrary auxiliary prime $\nu$, and use an algorithm from Wiedemann \cite{wiedemann} with some small modifications and additions. The changes to Wiedemann's algorithm that we make serve two purposes: some of them result in speedups for our problem specifically, and others are needed to guarantee that we find all $q$-expansions. In this section we'll outline these changes.

\subsection{Shifting eigenvalues}\label{charpoly_shift_subsection}
To compute the characteristic polynomial of $T_2$, we first compute the characteristic polynomial of $T_2 + kI$ for some integer $k$, and then make a change of variables. We do this for two reasons: to try and avoid singular matrices, and to give another parameter to modify in the random algorithm.\\
\\
Wiedemann's algorithm as described in \cite{wiedemann} is significantly more involved for singular matrices. In our case, if there's a newform in $S_2(p)$ with $a_2 = 0$, then the matrix representation of $T_2$ will be singular over $\Z$. To avoid these we instead work with the matrix $T_2 + kI$ for $k > 3$. The matrix $T_2$ over $\Z$ has one eigenvalue of $3$, and the others are guaranteed to be real and at most $2\sqrt{2}$ in absolute value. Thus $T_2 + kI$ is guaranteed to be nonsingular over $\Z$ whenever $k > 3$. The reduction of $T_2 + kI$ modulo $\nu$ might end up being singular anyway. We discuss this in section \ref{charpoly_variations_subsection}.

\subsection{Varying parameters}\label{charpoly_variations_subsection}
Wiedemann's algorithm is a random algorithm, and, for any given random input, fails a non-negligible amount of the time. Our purposes give us the freedom to vary two parameters which would be fixed in other situations: the shift $k$ (see section \ref{charpoly_shift_subsection}) and the modulus $\nu$. We got significant speedups by tweaking our algorithm to vary these parameters while also varying the random inputs.\\
\\
Overall, our implementation of Wiedemann's algorithm takes $4$ inputs:
\begin{itemize}
\item a random starting vector $u$,
\item a random coordinate $i$,
\item a shift $k$, and
\item a modulus $\nu$.
\end{itemize}
\noindent
To choose a random starting vector, we take the zero vector and set $50$ random entries to $1$. We found that setting only one entry to $1$ causes the algorithm to fail more often.\\
\\
The random coordinate $i$ which we use for the Berlekamp-Massey part of the algorithm is chosen uniformly at random.\\
\\
Even though the matrix $T_2 + kI$ is guaranteed to be nonsingular over $\Z$ whenever $k > 3$, the reduction mod $\nu$ might happen to be singular, essentially ``by chance''. In this case, the Berlekamp-Massey part of the algorithm will be given a power series which is not invertible. When this happens, we increment the shift $k$ in addition to choosing new values of $v$ and $i$.\\
\\
The cases which are the most computationally intensive, by a wide margin, are cases where the characteristic polynomial of $T_2$ has repeated factors. The runtime of our implementation of Wiedemann's algorithm is proportional to the multiplicity of the most frequently occurring factor (see section \ref{charpoly_minpolytocharapoly_subsection}), and, if the repeated factors are ones which might be reductions of factors of the characteristic polynomial over $\Z$, then we'll have to run the very expensive part of the algorithm which attempts to lift eigenspaces of dimension larger than $1$ (see section \ref{lift_section}). Genuine repeated factors of the characteristic polynomial $\Z$ are quite rare. Thus, we've found that the fastest approach is to run Wiedemann's algorithm for at most $2$ random $u$'s, and, if the algorithm fails for both choices, to change our small prime $\nu$ in case the failure was caused by a spurious repeated factor. After changing $\nu$, we increase the maximum number of choices of $u$, so that the algorithm does eventually find genuine high dimensional eigenspaces. Every time we change $\nu$, we have to recompute all of our iterates, so this results in a significant slowdown for levels which do have high dimensional eigenspaces, but there are very few of these.

\subsection{Getting the characteristic polynomial from the minimal polynomial}\label{charpoly_minpolytocharapoly_subsection}
For our purposes, it's important to find not just the minimal polynomial of $T_2$ (which is what Wiedemann's algorithm as described in \cite{wiedemann} yields), but the full characteristic polynomial. We need to do this for two reasons.\\
\\
First, we know the degree of the characteristic polynomial (because we know the dimension of $S_2(p)$), but not the degree of the minimal polynomial. Wiedemann's algorithm has some chance to fail to find the minimal polynomial, and when it does so, it outputs a polynomial which properly divides the minimal polynomial, and does not detect that it has failed to find the full minimal polynomial. We want to provably find all Galois orbits of $S_2(p)$, so it's necessary for us to be able to detect when the algorithm fails. By instead computing the characteristic polynomial, we can guarantee that the algorithm has succeeded by checking the degree.\\
\\
Second, a repeated factor of the characteristic polynomial of $T_2$ over $\Z$ indicates the presence of multiple Galois orbits of newforms with the same $a_2$. If we were to only compute the minimal polynomial of $T_2$, then there is some chance we would compute the $q$-expansion of only some of these newforms. In practice, this would be unlikely but not impossible; we elaborate on this in more detail in section \ref{lift_section}.\\
\\
Given the minimal polynomial of $T_2$, we use two techniques to produce the characteristic polynomial: comparing with known top coefficients (\ref{charpoly_topcoeff_subsubsection}), and checking eigenspace dimensions (\ref{charpoly_dimcheck_subsubsection}).

\subsubsection{Comparing with known top coefficients}\label{charpoly_topcoeff_subsubsection}
The coefficient of the second-highest degree term of the characteristic polynomial of a matrix $M$ is equal to $-\text{tr}(M)$. Similarly, the coefficient of the third-highest degree term of the characteristic polynomial is given by the expression \cite{mse}
\begin{equation}\label{thirdcoeffeqn}
  \sum_{1 \leq i < j \leq \text{dim}(M)} M_{ii}M_{jj} - M_{ij}M_{ji} .
\end{equation}
Computing the trace of any matrix $M$ takes a time of only $\cO(\text{dim}(M)^{1+\eps})$, and evaluating the expression (\ref{thirdcoeffeqn}) can be done in time $\cO(\text{dim}(M)^{2+\eps})$. Furthermore, because the matrix $T_2$ is sparse and nearly symmetric, these quantities can be computed more quickly. Expressions similar to the trace and (\ref{thirdcoeffeqn}) exist for other coefficients as well, and, while these would take too long to compute for general matrices, may be efficiently computable for sparse symmetric matrices like $T_2$. We haven't investigated this, but it would likely lead to a small improvement in the running time of our algorithm.\\
\\
If the characteristic polynomial is of degree at most $2$ more than than the degree of the minimal polynomial, and if the coefficients of the $3$ leading terms of both polynomials are known, then the ratio of these two polynomials can be found from an elementary calculation. We do this to find the characteristic polynomial whenever possible. Moreover, note that this trick still works even if the missing factors don't divide the minimal polynomial, so it also helps in some cases where Wiedemann's algorithm fails to find the minimal polynomial.

\subsubsection{Checking eigenspace dimensions}\label{charpoly_dimcheck_subsubsection}
This section of the overall algorithm is used only after Wiedemann's algorithm has been tried for $3$ or more initial random starting vectors $u$, which we'll label $u_1, \dots, u_n$. We stored the values of $T_2^ku_i$ for $0 \leq k \leq \text{dim}(T_2)$ and $1 \leq i \leq n$ (or, more precisely, the first $1000$ entries of these vectors; see section \ref{charpoly_iterates_subsection}). Let $\mu(t)$ denote the highest degree polynomial that's been returned by Wiedemann's algorithm thus far (so $\mu$ is a candidate for the minimal polynomial of $T_2$). Finding all the roots of $\mu$ takes time $\cO(p^{1+\eps})$ using the Fast Fourier Transform (henceforth ``FFT'') \cite{cooley-tukey}. Then, given a root $\lambda$ of $\mu$, the vector
$$v_{i,\lambda} \coloneqq \frac{\mu(T_2)}{T_2 - \lambda}u_i$$
is an eigenvector of $T_2$ with eigenvalue $\lambda$. Given the iterates of $u_i$, computing the first $1000$ entries of $v_{i,\lambda}$ takes time $\cO(p^{1+\eps})$. The dimension of the span of the vectors $v_{i,\lambda}$ is at most the dimension of the $\lambda$-eigenspace (and it's very likely that these dimensions will be equal, provided the number of iterates $n$ is at least the dimension of the eigenspace). This allows us to give lower bounds for the multiplicity with which linear factors of the minimal polynomial occur in the characteristic polynomial. In practice, the additional linear factors found with this trick, in combination with the trick from section \ref{charpoly_topcoeff_subsubsection}, are usually enough to determine the characteristic polynomial of $T_2$.

\subsection{Storing iterates}\label{charpoly_iterates_subsection}
To implement our version of Wiedemann's algorithm, we needed information about the iterates $T_2^ku$ of the random starting vector $u$ for three different purposes.
\begin{enumerate}
\item We computed these iterates as $T_2^ku = T_2(T_2^{k-1}u)$, so every coordinate of the previous iterate is stored temporarily.
\item To run the Berlekamp-Massey part of Wiedemann's algorithm, we used the $i^{\text{th}}$ coordinate of $T_2^ku$ for all $k < 2\text{dim}(T_2) + 10$, for an arbitrary $i$ (which we chose uniformly at random).
\item To implement the trick from section \ref{charpoly_dimcheck_subsubsection}, we used the first $1000$ coordinates of $T_2^ku$ for $k < \text{dim}(T_2)$. Here the choice of $1000$ is largely arbitrary; we just need to take enough coordinates to avoid spurious linear dependencies.
\end{enumerate}
The runtime of the entire algorithm depended heavily on computing these iterates quickly. Asymptotically, computing characteristic polynomials is the dominant term of the overall time complexity, and computing the iterates $T_2^ku$ was the step in our implementation of Wiedemann's algorithm that took the longest. Thus, we wanted to compute and store as few of the iterates as possible. Moreover, we structured our code so that python and Sage governed the algorithm at the top level, while computationally intensive sections were run in faster languages; these iterates were computed using cython directly. Because our code interfaced different programming languages, we needed to use the cython data to generate corresponding python data. If we store too much data, then this translation can take a significant amount of time. In the extreme case of storing the entire iterates $T_2^ku$, the translation took much longer than the actual computation of the data in cython.\\
\\
The other important reason to be judicious in how much of the iterates are stored is concern for the space complexity of the algorithm. We wanted our algorithm never use more than $\cO(p^{1+\eps})$ memory, and to accomplish this we were only able to store a constant number of coordinates per iterate.

\section{Finding $\Z$-eigenbases}\label{lift_section}
Let $\chi_\Z$ be the characteristic polynomial of $T_2$ over $\Z$, and $\chi_\nu$ its reduction modulo $\nu$. At this point, the algorithm has computed $\chi_\nu$ for one choice of $\nu$. The objective of this section is to determine all irreducible factors of degree $6$ or less of $\chi_\Z$, and, for each factor, an eigenbasis defined over $\Z$ and a simultaneous eigenvector (defined over a number field) for all Hecke operators. Here we use the term ``eigenspace'' and related terms to mean $\text{ker}(\rho(T_2))$ for some irreducible factor $\rho$ of $\chi_\Z$. We'll need $\Z$-bases of these spaces to use a formula from Mestre's M\'{e}thode des Graphes, which we'll discuss in section \ref{qexpansion_section}.\\
\\
If an irreducible factor $\rho$ of $\chi_\Z$ divides $\chi_\Z$ exactly once, then $\rho$ corresponds to a single Galois orbit of a newform of $S_2(p)^\pm$, where $\pm$ denotes the Atkin-Lehner eigenvalue. If an irreducible factor $\rho$ divides $\chi_\Z$ twice or more then the situation is more complicated. Multiple Galois orbits of newforms with the same $a_2$ lead to repeated factors of $\chi_\Z$, as do newforms whose Hecke field is of strictly larger degree than the degree of $a_2$. In these cases, Galois orbits of newforms will correspond to the minimal nontrivial subspaces of $\text{ker}(\rho(T_2))$ which are invariant under the action of every $T_\ell$ simultaneously.\\
\\
The key idea behind the lifting algorithm described in this section is taking advantage of the heuristic that the coordinates in the vectors making up the $\Z$-eigenbasis are very likely to be small in absolute value. We pick some number of ``candidate lifts'' of an $\F_\nu$-eigenbasis in which the most common entries are small integers, and then check directly whether or not our candidate lifts are eigenvectors of $T_2$ over $\Z$. In a large majority of cases, if $\rho$ does indeed divide $\chi_\Z$, then checking only a handful of candidate lifts is enough to find a $\Z$-eigenbasis.

\subsection{Finding eigenvalues}\label{eigenvalues_subsection}
Let $\{\rho_i\}$ denote the set of polynomials which occur as the minimal polynomial of $a_2$ of a simple Abelian variety over $\F_2$ of dimension $6$ or less. There are $96795$ such polynomials, and they can be found in the LMFDB. Every irreducible factor of $\chi_\Z$ of degree $6$ or less is necessarily one of these.\\
\\
To find the irreducible factors of $\chi_\Z$, we first determine which $\rho_i$'s divide $\chi_\nu$. We do this by iterating over the set $\{\rho_i\}$ and checking for divisibility one by one, but one could use FFT \cite{cooley-tukey} to do this more quickly if necessary.\\
\\
It is possible for $\rho_i$ to divide $\chi_\nu$ but not $\chi_\Z$. When we find a $\rho_i$ which divides $\chi_\nu$, we attempt to produce an eigenbasis defined over $\Z$ for this factor using a method which we describe in the rest of this section. If we succeed, this proves that $\rho_i$ does divide $\chi_\Z$. Our method for producing a $\Z$-eigenspace is not guaranteed to work, though it almost always does in practice. In the cases where it doesn't, we compute $\chi_{\nu'}$ for some other small prime $\nu' \neq \nu$, and see if $\rho_i$ divides $\chi_{\nu'}$. If $\rho_i\nmid\chi_{\nu'}$, then we've proven that $\rho_i\nmid\chi_\Z$, and, conversely, if $\rho_i\nmid\chi_\Z$, then there is guaranteed to be some prime $\nu'$ for which $\rho_i\nmid\chi_{\nu'}$ (and usually only one additional prime $\nu'$ needs to be checked).\\
\\
It would be difficult to determine the Galois orbits of size $7$ or greater using this approach, because the number of simple Abelian varieties over $\F_2$ of given dimension grows very quickly. In section \ref{highdeg_section}, we describe the algorithm we used to determine the size of every Galois orbit, including the ones of size $7$ or more. For every prime level between $10,\!000$ and $1,\!000,\!000$ and each Atkin-Lehner eigenspace, the characteristic polynomial $\chi_\Z$ had only one irreducible factor of degree $7$ or more. For prime levels between $1,\!000,\!000$ and $2,\!000,\!000$ we computed the $q$-expansions of the newforms of degree $6$ or less but did not investigate the decomposition of the rest of the space.

\subsection{Lifting $1$-dimensional eigenspaces}\label{lift_1dim_subsection}
When an eigenspace is $1$-dimensional (which is the case exactly when $\rho(t) = t - \lambda$ is linear and divides $\chi_\nu$ exactly once), our algorithm is fairly straightforward. First, we find an eigenvector $v$ of $T_2$ by computing
$$v \coloneqq \frac{\mu(T_2)}{T_2 - \lambda}u$$
for some random starting vector $u$, where $\mu$ is the minimal polynomial of $T_2$ on the given Atkin-Lehner eigenspace. Computing $v$ takes time $\cO(p^{2+\eps})$, so, for the levels for which require that we attempt to lift eigenspaces, this is a leading term of the over asymptotic time complexity of our algorithm, and these levels presumably make up a small but strictly positive proportion of all levels. The computation of $v$ in this way involves successively computing the iterates $T_2^ku = T_2(T_2^{k-1}u)$. This computation was done previously in the section of our algorithm which used Wiedemann's algorithm, but it was impossible to store these iterates while meeting our goal of having a space complexity of $\cO(p^{1+\eps})$. Implementations of our algorithm which can afford to store the iterates from the Wiedemann section can use those iterates here, saving some time.\\
\\
With our eigenvector $v$ defined over $\F_\nu$, we generate our list of ``candidate lifts'' of $v$ to an eigenvector $\hat{v}$ defined over $\Z$ by guessing that the most likely scenario is that the most common nonzero entry of $\hat{v}$ is $1$, the second most likely scenario is that the most common entry is $2$, and so on. Note that it suffices to check for strictly positive most common entries, since both $\hat{v}$ and $-\hat{v}$ span the $\Z$-eigenbasis.\\
\\
Let $\alpha \in \F_\nu$ denote the most common nonzero entry of $v$. Then, our candidate lifts are $\frac{1}{\alpha}v, \frac{2}{\alpha}v, \dots$, where each $\F_\nu$ coordinate of these vectors is lifted to the integer of smallest absolute value in that residue class. As these candidate lifts are generated, we multiply them by $T_2$ to check directly whether or not they're eigenvectors, and return the first candidate lift that is an eigenvector. This lift is guaranteed to have entries with a $\text{gcd}$ of $1$, so it spans the $\Z$-eigenspace (a condition which is necessary for the part of the algorithm described in section \ref{qexpansion_section} to succeed). It is also guaranteed to be a simultaneous eigenvector of all the Hecke operators.\\
\\
Our implementation of this algorithm generates $50$ candidate lifts, and then, if no lift to $\Z$ has been found, gives up and declares it has failed to lift the $\F_\nu$ eigenvector. This prompts the algorithm to try a different small prime $\nu' \neq \nu$, as discussed in section \ref{eigenvalues_subsection}. We encountered no cases in which a lift existed but was not found.

\subsection{Lifting higher dimensional eigenspaces}\label{lift_highdim_subsection}

\subsubsection{Finding an $\F_\nu$-eigenbasis}\label{lift_highdim_local_subsubsection}
We begin like we did in the $1$-dimensional case, by computing
$$v \coloneqq \frac{\mu(T_2)}{\rho(T_2)}u$$
for some random starting vector $u$, where $\mu$ is the minimal polynomial of $T_2$ on the given Atkin-Lehner eigenspace. The vector $v$ is in the $\F_\nu$-kernel of $\rho(T_2)$, and is what we're calling an eigenvector (because it's an eigenvector up to Galois conjugacy for the field generated by $\rho$). The same discussion as that in section \ref{lift_1dim_subsection} applies here: computing this eigenvector is a small but nonzero part of the leading term in asymptotic time complexity, and the iterates from the Wiedemann component of our algorithm can't be used to save time in the computation in our implementation because of space limitations.\\
\\
Let $r$ denote the multiplicity with which $\rho$ divides $\chi_\nu$. If $r = 1$, then we compute the dimension of the span of the vectors
$$v, T_2v, T_2^2v, \dots, T_2^{\text{deg}(\rho)-1}v.$$
If the dimension of this span is equal to $\text{deg}(\rho)$, then these vectors form an $\F_\nu$-eigenbasis of the eigenspace. If this dimension is strictly less than $\text{deg}(\rho)$, or if $r>1$, then, starting from $\ell=3$, we compute matrix which corresponds to the action of $T_\ell$ on the given Atkin-Lehner eigenspace by exploring the supersingular $\ell$-isogeny graph in the way described in section \ref{T2_section}, and compute the dimension of the span of the vectors
$$v, T_\ell v, T_\ell^2v, \dots, T_\ell^{r\text{deg}(\rho)-1}v,$$
stopping our iteration over $\ell$ when this dimension is equal to $r\deg(\rho)$. We also periodically recompute $v$ with a new random choice of $u$, because it's possible (but unlikely) that $v$ happens to be a simultaneous eigenvector of all Hecke operators, which would cause our iteration over $\ell$ to never terminate (in principle one could add a clause in our algorithm to take advantage of this whenever it happened, but in practice this never happened). In doing this, we find an $\ell$ such that all newforms of $S_2(p)^\pm$ have distinct $a_\ell$'s. We need these $a_\ell$'s to be distinct later, essentially so that we can separate the $\Z$-eigenbases. Let $M$ denote the $\text{dim}(S_2(p)^\pm) \times r\text{deg}(\rho)$ matrix whose columns are $v, T_\ell v, T_\ell^2v, \dots, T_\ell^{r\text{deg}(\rho)-1}v$.

\subsubsection{Finding a $\Z$-eigenbasis for the full $\rho$-eigenspace}\label{lift_highdim_fullspace_subsubsection}
We now make use of our key heuristic that our $\Z$-eigenbases are very likely to be made up of vectors which have small entries. To do this, we start by making two lists (which we generate as needed, as opposed to storing in memory):
\begin{itemize}
\item $\cM$: a list of the most common rows of $M$ with the condition that any rows we add to this list are not in the span of the rows we've added previously, and
\item $\cC$: a list of ``candidate columns'', which are $r\text{deg}(\rho)$-tuples with small integer entries, ordered in such a way such that, generally speaking, tuples with smaller entries are listed first.
\end{itemize}
Our implementation requires us to pick an ordering of the set $\cM^{r\text{deg}(\rho)} \times \cC$ so that we may iterate over it. We do this more or less arbitrarily. Experimentally, we found that assigning a ``size'' to elements of $\cC$ which was proportional to the sum of squares of the entries, assigning a ``size'' to elements of $\cM^{r\text{deg}(\rho)}$ which was proportional to the sum of the squares of the inverse appearance frequencies, and then ordering $\cM^{r\text{deg}(\rho)} \times \cC$ by the product of these sizes lead to having to check fewer candidate lifts than lexicographic orderings.\\
\\
We iterate over $\cM^{r\text{deg}(\rho)} \times \cC$. For each element $(m,c)$, let $L$ be the linear combination of columns of $M$ which yields $c$ when restricted to the rows in $m$. We then produce a ``candidate lift'' $\hat{v}$, which is the vector with integer entries produced by taking the combination $L$ of the columns of $M$, and then lifting each coordinate of the resulting $\F_\nu$ vector to the integer in the appropriate residue class with smallest absolute value. We then check directly whether or not $\hat{v}$ is in $\text{ker}(\rho(T_2))$. We continue this process until we've found $r\text{deg}(\rho)$ linearly independent candidate lifts. This set of candidate lifts forms a $\Z$-basis of $\text{ker}(\rho(T_2))$, but our algorithm requires a $\Z$-basis for each Galois orbit, so it remains to decompose our $\Z$-basis in this way. We do this next.

\subsubsection{Finding a $\Z$-eigenbasis for each Galois orbit}\label{lift_highdim_orbit_subsubsection}
The method we describe in this section for finding a $\Z$-eigenbasis for each Galois orbit is built around the observation that the action of $T_\ell$ fixes $\text{ker}(\rho(T_2))$. This means that each column of $T_\ell M$ can be expressed as a linear combination of the columns of $M$. Let $S$ denote the $r\text{deg}(\rho) \times r\text{deg}(\rho)$ integer matrix whose entries are the coefficients of these linear combinations (which are all small, since $T_\ell$ is given by a matrix with at most $2(\ell+1)$ nonzero entries per row, each of which is $\pm 1$).\\
\\
Let $\chi_S$ denote the characteristic polynomial of $S$, and let $h_1, h_2, \dots$ denote the minimal polynomials of the $a_\ell$'s of the newforms we're considering. As discussed in section \ref{lift_highdim_local_subsubsection}, we chose $\ell$ in a way that guarantees that each of these minimal polynomials is distinct. By construction, we have
$$\chi_S = \prod_i h_i.$$
We can then obtain, with some elementary linear algebra, a $\Z$-basis for $\text{ker}(h_i(S))$ and a set of eigenvectors of $S$ (which will be defined over number fields). Taking the linear combinations of the columns of $M$ whose coefficients are given by the elements of these $\Z$-bases and eigenvectors then yields the $\Z$-eigenbases and simultaneous eigenvectors we need. Our implementation requires only one simultaneous eigenvector per Galois orbit, which we pick in an arbitrary way.

\section{Computing $q$-expansions}\label{qexpansion_section}
Let $v$ be a simultaneous eigenvector of all the Hecke operators whose coordinates are indexed by the supersingular $j$-invariants over $\overline{\F}_p$. In \cite{mestre}, Mestre gives the $q$-expansion of the associated newform $f$ modulo any prime $\mathfrak{p}$ above $p$ in $K \coloneqq \Q(a_2, a_3, \dots)$ via an equality of power series involving the $q$-expansion of the modular $j$-function $j(q)$:
\begin{equation}\label{mestre_eqn}
  \left(\sum_j v_jj\right)f(q)\frac{dq}{q} = \sum_j v_j\frac{dj(q)}{j(q) - j} \quad\text{ mod }\mathfrak{p}.
\end{equation}
The Weil bound then allows one to determine $a_\ell$ for all $\ell \ll p^2$. Our algorithm uses this formula to compute $a_n$ for $n$ up to the Sturm bound.\\
\\
Section \ref{qexpansion_zbasis_subsection} gives a version of equation (\ref{mestre_eqn}) that does computations over $\F_p$ using the $\Z$-eigenbases computed earlier. Sections \ref{qexpansion_jq_subsection} and \ref{qexpansion_powerseries_subsection} detail how we evaluate the right hand side of (\ref{mestre_eqn}) efficiently.

\subsection{$\Z$-eigenbasis version of Mestre's identity}\label{qexpansion_zbasis_subsection}
We use the simultaneous eigenvector $v$ to compute the values $a_{\ell_1}, a_{\ell_2}, \dots, a_{\ell_{\text{deg}K}}$ (which might require computing new Hecke matrices using supersingular isogeny graphs). With each vector $u_k$ in the $\Z$-eigenbasis, we evaluate the right hand side of (\ref{mestre_eqn}), and let $\psi_k(q)$ denote the resulting power series. As we'll discuss in section \ref{qexpansion_powerseries_subsection}, the power series $\psi_k(q)$ is guaranteed to be defined over $\F_p$. We write $\psi_k[\ell]$ to refer to the coefficient of $q^\ell$ in $\psi_k(q)$. Pick a basis $\{r_i\}$ of the ring of integers $\cO_K$. Write each $a_\ell$ in terms of this basis:
$$a_\ell \eqqcolon \sum_i\alpha_i[\ell] r_i.$$
Then, from (\ref{mestre_eqn}), it follows that there exist coefficients $\beta_{i,k} \in \F_p$ such that, for every $\ell$, we have the a linear combination
\begin{equation*}
  a_\ell = \sum_{i,k}\beta_{i,k}\psi_k[\ell]r_i \quad\text{mod }\mathfrak{p},
\end{equation*}
so
\begin{equation}\label{mestre_matrix_eqn}
  \begin{pmatrix}
    \alpha_1[\ell_1] & \alpha_1[\ell_2] & \dots\\
    \alpha_2[\ell_1] & \alpha_2[\ell_2] & \dots\\
    \vdots & \vdots & \ddots
  \end{pmatrix}
  =
  \begin{pmatrix}
    \beta_{1,1} & \beta_{1,2} & \dots\\
    \beta_{2,1} & \beta_{2,2} & \dots\\
    \vdots & \vdots & \ddots
  \end{pmatrix}
  \begin{pmatrix}
    \psi_1[\ell_1] & \psi_1[\ell_2] & \dots\\
    \psi_2[\ell_1] & \psi_2[\ell_2] & \dots\\
    \vdots & \vdots & \ddots
  \end{pmatrix}
  \quad\text{mod }p.
\end{equation}
We solve the matrix equation (\ref{mestre_matrix_eqn}) for the coefficients $\beta_{i,k}$ (as long as all the matrices involved are invertible; compute $a_{\ell}$ for more $\ell$'s if there happens to be a linear dependence). We can then find the $q$-expansion of our newform as
$$f(q)\frac{dq}{q} = \sum_i\left(\sum_k\beta_{i,k}\psi_k[\ell]\right)r_i q^\ell.$$
We lift the coefficient $\sum_k\beta_{i,k}\psi_k[\ell]$ to $\Z$ by ensuring that the Hecke eigenvalues satisfy the Weil bound. Usually the coefficient lifts to the integer in the residue class of smallest absolute value if one chooses $\{r_i\}$ to be a reduced basis of $\cO_K$.

\subsection{Computing $j(q) \text{ mod }p$}\label{qexpansion_jq_subsection}
To evaluate the right hand side of (\ref{mestre_eqn}) up to the Sturm bound, we first need to compute $\cO(p)$ coefficients of $j(q)$ modulo $p$. We compute these coefficients using the identity
\begin{equation}\label{j_eta_eqn}
  j(q) = \frac{E_{12}(q)}{\left(\eta(q)^3\right)^8} - \frac{82104}{691} + 744,
\end{equation}
where $E_{12}(q)$ is the weight $12$ classical Eisenstein series (normalized to have constant term $1$) and $\eta(q)$ is the Dedekind $\eta$-function. We chose equation (\ref{j_eta_eqn}) because the function $\eta(q)^3$ has a ``sparse'' $q$-expansion:
$$\eta(q)^3 = \sum_{k=0}^\infty (-1)^k (2k+1) q^{\frac{k(k+1)}{2}}.$$
Computing the first $p$ terms of the quotient of two power series takes time $\cO(p^{1+\eps})$ using FFT, and empirically dividing $8$ times by $\eta(q)^3$ was faster than alternatives, such as dividing by $\Delta(q)$ or $24$ times by $\eta(q)$.\\
\\
In section \ref{qexpansion_powerseries_subsection} we'll also need the power series $j'(q)$. Differentiating $j(q)$ takes time $\cO(p^{1+\eps})$.

\subsection{Fast power series algorithms}\label{qexpansion_powerseries_subsection}
\subsubsection{Atkin-Lehner eigenspace to Mestre's formula}\label{qexpansion_unfolding_subsubsection}
The eigenvectors $\{u_k\}$ we computed are indexed by pairs $(j,j^\sigma)$ with $j < j^\sigma$ (in the $+$ space) or $j \leq j^\sigma$ (in the $-$ space), where $<$ is the arbitrary but fixed ordering from section \ref{T2_section}. To use Mestre's formula (\ref{mestre_eqn}), we need to construct eigenvectors $\{v_k\}$ that are indexed by supersingular $j$-invariants. Recall that in section \ref{T2_section} we chose a generator $\xi$ of $\F_{p^2}$ which satisfied $\xi^\sigma = -\xi$. This allows us to construct each $v_k$ easily: if $j \not\in \F_p$ we set $v_k(j) = u_k((j,j^\sigma))$ and $v_k(j^\sigma) = \mp u_k((j,j^\sigma))$, and if $j \in \F_p$ we set $v_k(j) = (1 \mp 1)u_k((j,j^\sigma))$. Then $v_k$ is defined over $\Z$ (because $u_k$ is), and the sum in (\ref{mestre_eqn}) is either Galois-invariant or Galois anti-invariant. In the latter case we then divide by $\xi$.
  
\subsubsection{Evaluating Mestre's formula}\label{qexpansion_binarytree_subsubsection}
Evaluating the sum of rational functions of $q$-expansions in $(\ref{mestre_eqn})$ by evaluating each term separately and then then adding would take time $\gg p^2$. However, it is possible to evaluate this sum in time $\cO(p^{\frac{3}{2}+\eps})$ by doing a ``binary tree decomposition'' and using a power series composition algorithm from Brent and Kung \cite{brent-kung}.\\
\\
Given constants $\gamma_1, \dots, \gamma_M$ and $j_1, \dots, j_M$, let $P$ and $Q$ be polynomials such that
$$\frac{P(x)}{Q(x)} \coloneqq \sum_{i=1}^M \frac{\gamma_i}{x - j_i}.$$
Define
$$S(x,y) \coloneqq \sum_{i=x}^y \frac{\gamma_i}{x - j_i}.$$
We compute $P$ and $Q$ recursively as
\begin{align*}
  \sum_{i=1}^M \frac{\gamma_i}{x - j_i} &= S(1,M)\\
  &= S\!\left(1,\tfrac{M}{2}\right) + S\!\left(\tfrac{M}{2}+1,M\right)\\
  &= \left[S\!\left(1,\tfrac{M}{4}\right) + S\!\left(\tfrac{M}{4}+1,\tfrac{M}{2}\right)\right] + \left[S\!\left(\tfrac{M}{2}+1,\tfrac{3M}{4}\right) + S\!\left(\tfrac{3M}{4}+1,M\right)\right]\\
  &= \dots\\
  &= \sum_{i=1}^M S(i,i).
\end{align*}
We start from the bottom expression and work upwards. At the $k^{\text{th}}$ step we compute $\cO\!\left(\frac{M}{2^k}\right)$ sums of two terms. Each term at the $k^{\text{th}}$ step is a rational function whose numerator and denominator have degree $\cO(2^k)$. Multiplying two polynomials of degree $d$ takes time $\cO(d^{1+\eps})$ using FFT. Thus, each step takes time $\cO\!\left((2^k)^{1+\eps}\frac{M}{2^k}\right) = \cO(M^{1+\eps})$. There are $\cO(\log M)$ steps in this procedure. Thus, computing $P$ and $Q$ takes time $\cO(M^{1+\eps}\log M) = \cO(M^{1+\eps})$.\\
\\
In the context of our problem, this means that we can, in time $\cO(p^{1+\eps})$, compute polynomials $P, Q$ of degree $\cO(p)$ such that
$$\sum_j v_j \frac{dj(q)}{j(q) - j} = \frac{P(j(q))}{Q(j(q))}j'(q)dq.$$
Writing $\frac{1}{Q(j(q))}$ as a power series in $j(q)$ with $\cO(p)$ terms of precision takes time $\cO(p^{1+\eps})$, and computing the product $R(j(q)) \coloneqq P(j(q))\frac{1}{Q(j(q))}$ as a power series in $j(q)$ with $\cO(p)$ terms of precision also takes time $\cO(p^{1+\eps})$.\\
\\
Using \cite{brent-kung}, we compute the composition $R(j(q))$ to $\cO(p)$ terms of precision in time $\cO(p^{\frac{3}{2}+\eps})$. Finally, multiplying by $j'(q)dq$ takes time $\cO(p^{1+\eps})$.

\section{Checking for high degree factors of the characteristic polynomial}\label{highdeg_section}
As mentioned in section \ref{eigenvalues_subsection}, our algorithm doesn't find newforms whose $a_2$ has minimal polynomial of degree $7$ or more. It's of value to at least determine the dimensions of all of the newforms of $S_2(p)$, even if we can't compute all of their $q$-expansions. This part of the algorithm is by far the most expensive time and space wise, and is the most technically involved. This part of the algorithm isn't required for computing $q$-expansions and can be omitted.\\
\\
Our algorithm is designed to try and efficiently determine that there is exactly one irreducible factor of $\chi_\Z$ of degree $7$ or more, since this was the case for every prime level between $10,\!000$ and $1,\!000,\!000$. The approach we take is based on the straightforward observation that if, for some small prime $\nu$, the characteristic polynomial $\chi_\nu$ has no factor of degree $d$, then $\chi_\Z$ cannot have a factor of degree $d$ either. Thus, for each possible degree $d$ between $7$ and $\thalf \text{deg}\chi_\Z$ we aim to find some small prime $\nu$ such that $\chi_\nu$ has no factor of degree $d$.

\subsection{Factoring the characteristic polynomial modulo many small primes}\label{highdeg_factor_subsection}
The first part of the algorithm will require finding and factoring characteristic polynomials $\chi_{\nu_1}, \chi_{\nu_2}, \dots$ of $T_2$ on a given Atkin-Lehner eigenspace modulo many small primes, which we'll do as needed in what follows. Our algorithm is slightly more efficient if we replace $\chi_{\nu_i}$ with $\chi_{\nu_i}$ divided by all of the irreducible factors of degree $6$ or less which we know divide $\chi_\Z$. Factoring polynomials over finite fields has high time complexity in theory, but experimentally we found that factoring $\chi_{\nu_i}$ appeared to take time $\cO(p^{2+\eps})$ in our specific case. This is still the most significant chunk of the runtime, but it's not problematically expensive.\\
\\
The reason factoring is feasible can be explained heuristically. There exists an algorithm based on FFT which factors polynomials over finite fields in three steps \cite{cohen}:
\begin{enumerate}
\item Eliminating square factors of $\chi_{\nu_i}$
\item Factoring $\chi_{\nu_i}$ into products of irreducible polynomials of equal degree
\item Factoring each of these products of irreducibles of equal degree
\end{enumerate}
The first step is not problematic because we obtain the squarefree part of $\chi_{\nu_i}$ (which is the minimal polynomial) directly through Wiedemann's algorithm. Even ignoring this, eliminating the square factors of $\chi_{\nu_i}$ can be done efficiently by computing $\text{gcd}(\chi_{\nu_i}, \chi_{\nu_i}')$.\\
\\
The second step can be done in time $\cO(p^{2+\eps})$ using FFT.\\
\\
The third step is the one which is theoretically challenging, but, in practice, it's very rare that $\chi_{\nu_i}$ has multiple factors of the same degree if that degree is large. Thus, we can use Rabin's irreducibility test \cite{rabin}, which runs in time $\cO(p^{2+\eps})$, to determine which of the products from the second step require factoring. The products that do require factoring end up being products of a handful of small degree factors which can be factored quickly. There is also never a requirement to use any specific $\nu$, so at worst one could abandon computations with $\nu$'s that were stuck on this step.\\
\\
With knowledge of this factorization we can use the algorithm described in section \ref{highdeg_sieve_subsection}. We've found that, in our implementation, our algorithm ran more quickly if we only continued using this factorization if it didn't have ``too many'' irreducible factors, with a threshold determined experimentally and in a fairly arbitrary way.

\subsection{Sieving possible degrees}\label{highdeg_sieve_subsection}
Let $E_i$ denote the set subset of $[7,\frac{1}{2}\text{deg}(\chi_\Z)] \cap \Z$ which we have yet to show cannot be degrees of factors of $\chi_\Z$ after running this part of the algorithm for $\nu_i$.  For each $i$, we determine
\begin{enumerate}
\item which elements of $E_i$ occur as the degrees of (not necessarily irreducible) factors of $\chi_{\nu_i}$, and,
\item for each $d$ in $E_i$, if there are only ``a few'' factors of $\chi_{\nu_i}$ with this degree, we record them to later use in the part of the algorithm described in section \ref{highdeg_weil_subsection}. We specify what we mean by ``a few'' later in this section.
\end{enumerate}
Define
$$D_i \coloneqq \{d \in E_i \,:\, \text{there exists a factor of }\chi_{\nu_i}\text{ with degree }d\}.$$
Given a factorization of $\chi_{\nu_i}$ into irreducibles, computing the set $D_i$ is a manifestation of the well-studied subset sum problem. Direct enumeration of all possible combinations of irreducible factors of $\chi_{\nu_i}$ takes time and space exponential in the number of irreducible factors and was infeasible in practice. There's a dynamic programming algorithm \cite{clrs} which computes $D_i$ in time $\cO(p^2)$ and space $\cO(p)$, but says nothing about what the corresponding factors are. We'll modify this dynamic programming algorithm so that it records which products of irreducible factors have degree $d$ whenever there are only ``a few'' such products. Doing this allows us to use the the method described in section \ref{highdeg_weil_subsection}, which leads to having to compute fewer characteristic polynomials $\chi_{\nu_i}$ and their factorizations. If one chooses, one can only compute $D_i$ and never anything about the factors themselves, and then use the dynamic programming algorithm directly. We've found that omitting the method from \ref{highdeg_weil_subsection} usually takes longer but requires less space for our implementation.\\
\\
To compute $D_i$, we create a variable $\Delta_i$ which will ultimately be the function with domain $D_i$ for which $\Delta_i(d)$ is either a set of combinations of irreducible factors whose product has degree $d$, or ``null'' if this set would have more than ``a few'' elements in it. We initialize $\Delta_i$ as the function with $\Delta_i(0) = \emptyset$ and no other elements in its domain.\\
\\
For $k = 1,2,\dots$, let $h_k$ denote the irreducible factors of $\chi_{\nu_i}$. We sort these factors in descending order according to their degrees (breaking ties arbitrarily), and compute the partial sums $P_k \coloneqq \sum_{\kappa\geq k} h_\kappa$. For each $k$, let $R_k$ denote the set of elements of $E_i$ which are not yet in $\Delta_i$. We iterate over either the domain of $\Delta_i$, or the set
$$\left\{d \in \Z \,:\, \text{min}(R_k)-P_k \leq d \leq \text{max}(R_k)\right\},$$
whichever is smaller. For each $d$, if $d \in \text{dom}(\Delta_i)$, then we check whether or not $d+\text{deg}(h_k)$ is in $\text{dom}(\Delta_i)$. If it isn't, then we set
$$\Delta_i(d+\text{deg}(h_k)) = \{H \cup \{h_k\}\,:\,H \in \Delta_i(d), h_k \not\in H\}.$$
If it is, then, if $\Delta_i(d+\text{deg}(h_k)) \neq \text{``null''}$, we compute the set
$$\cH' \coloneqq \Delta_i(d+\text{deg}(h_k)) \cup \{H \cup \{h_k\}\,:\,H \in \Delta_i(d), h_k \not\in H\}.$$
If $\cH'$ is larger than some fixed parameter $\eta$ given as input to the algorithm, we then replace $\cH'$ with ``null''. We took $\eta = 5$ in our implementation. Then we update $\Delta_i$ by setting $\Delta_i(d+\text{deg}(h_k)) = \cH'$.\\
\\
Previously, when we said we wouldn't find which products of irreducible factors multiplied together to give factors of the given degree if there were more than ``a few'' such products, this process of replacing values of $\Delta_i$ by ``null'' whenever the values would have cardinality more than $\eta$ is the condition we were referring to: at no intermediate step were there more than $\eta$ ways of obtaining that intermediate degree. It was necessary to have some bound of this sort, since otherwise we are directly enumerating all possible products of irreducible factors, which is infeasible. We chose $\eta = 5$ because empirically it was the best balance we found between yielding non-null values of $\Delta_i$ and not using excessive memory.\\
\\
After the iteration over $k$ finishes, we set $E_{i+1} = E_i \cap \text{dom}(\Delta_i)$. Moreover, for each $d \in E_{i+1}$, we store the values of $\Delta_i(d)$ whenever they're not ``null'', and then run the part of our algorithm described in section \ref{highdeg_weil_subsection} before continuing our iteration over $i$.

\subsection{Checking the Weil bound}\label{highdeg_weil_subsection}
If $\theta(t) = t^d + \theta_1t^{d-1} + \dots + \theta_d$ is a polynomial which divides $\chi_\Z$, then the coefficients of $\theta$ satisfy the Weil bound:
$$|\theta_j| \leq {d\choose j}(2\sqrt{2})^j.$$
This bound gives us a way to rule out the existence of a lift of a polynomial $\bar{\theta}(t) \in (\Z/m)[t]$ to a factor of $\chi_\Z$: if $\bar{\theta}$ has a coefficient which has no integer lifts that satisfy the Weil bound, then we know that $\bar{\theta}$ cannot lift to a factor of $\chi_\Z$. For our range of levels and our choice of small primes $\nu$, this Weil bound condition ends up being trivial for $j \geq 3$, but for $j = 1,2$ we use this bound to avoid computing and factoring more $\chi_{\nu_i}$'s than we would have to otherwise.\\
\\
For each degree $d \in E_i$, if there are primes $\nu_{i_1}, \nu_{i_2}, \dots$ for which
\begin{enumerate}
\item all ways to obtain a factor of degree $d$ from a product of irreducible factors of $\chi_{\nu_{i_k}}$ are known, and
\item $\prod_k \nu_{i_k}$ is large enough for the Weil bound strategy outlined above to not be trivial,
\end{enumerate}
then we use the Chinese Remainder Theorem on all possible combinations of products of irreducible factors to produce a list of ``candidate factors'' modulo $m = \prod_k \nu_{i_k}$. We then use the Weil bound as outlined above on these candidate factors one at a time to try and rule them out. Even if not all candidate factors can be ruled out, usually some can, and we remove these candidates so that this method is more likely to succeed as we continue to iterate over $i$.

\section{Results}
We used the algorithm described in this paper to compute the $q$-expansions of all weight $2$ eigenforms of dimension $g \leq 6$ and prime level between $10^4$ and $2\cdot 10^6$. The $q$-expansions of forms with level less than $10^4$ were computed by Best, Bober, Booker, Costa, Cremona, Derickx, Lowry-Duda, Lee, Roe, Sutherland, and Voight \cite{bbbccdldlrsv}. Forms with $g = 1$ correspond to elliptic curves, and Weierstrass equations for all elliptic curves of prime conductor less than $2\cdot 10^9$ were computed by Bennett, Gherga, and Rechnitzer \cite{bgr}. Both of these datasets are in the LMFDB \cite{lmfdb}.\\
\\
Below, we tabulate the number of forms of prime level in the ranges $[1,10^4]$, $[10^4, 10^6]$, and $[10^6, 2\cdot 10^6]$, grouped by the discriminant $\Delta$ of their Hecke fields and omitting discriminants which don't appear in our dataset.\\
~\\
\begin{tabular}{|c|r|r|r|r|l|}
  \multicolumn{2}{c}{} & \multicolumn{3}{c}{Number with prime level in...} & \multicolumn{1}{c}{} \\
  \hline
  \rule{0pt}{1em}$g$ & \multicolumn{1}{c|}{$\Delta$} & $[1, 10^4]$ & $[10^4, 10^6]$ & $[10^6, 2\cdot 10^6]$ & \multicolumn{1}{c|}{Levels} \\ [0.05ex]
  \hline
  \rule{0pt}{1em}\multirow{1}{*}{$1$}
  & $1$ & $329$ & $8843$ & $6406$ & $11\, \dots\, 1999957$ \\ [0.05ex]
  \hline
  \rule{0pt}{1em}\multirow{6}{*}{$2$}
  & $5$  & $158$ & $1900$ & $986$ & $23\, \dots\, 9973, 10103, 10267\, \dots\, 1997773, 1999867$ \\
  & $8$  & $37$  & $242$  & $100$ & $29\, \dots\, 9613, 12619, 13537\, \dots\, 1986043, 1991489$ \\
  & $13$ & $13$  & $40$   & $6$   & $73\, \dots\, 9967, 15193, 23473\, \dots\, 1773361, 1863347$ \\
  & $12$ & $1$   & $14$   & $3$   & $113, 13763, 15083, 15919, 22481\, \dots\, 1493441, 1894043$ \\
  & $21$ & $1$   & $3$    & $1$   & $1283, 112289, 329671, 577807, 1670563$ \\
  & $17$ &       & $1$    &       & $75653$ \\ [0.05ex]
  \hline
  \rule{0pt}{1em}\multirow{7}{*}{$3$}
  & $49$  & $34$ & $90$ & $30$ & $97\, \dots\, 9857, 12569, 13121\, \dots\, 1929943,1972423$ \\
  & $229$ & $8$  & $20$ & $1$  & $211\, \dots\, 6997, 11197, 14563\, \dots\, 614659, 1972651$ \\
  & $81$  & $3$  & $13$ &      & $127, 3581, 8513, 10753, 35591\, \dots\, 277793, 336551$ \\
  & $169$ & $2$  & $6$  & $3$  & $1481, 6569, 22943, 42209\, \dots\, 1221239, 1856201$ \\
  & $257$ & $9$  & $6$  & $1$  & $71\, \dots\, 26713, 39089, 224057, 255877, 463249, 581657, 1120969$ \\
  & $148$ & $12$ & $6$  &      & $41\, \dots\, 2341, 14929, 31039, 133117, 319489, 429397, 707801$ \\
  & $321$ & $2$  & $1$  &      & $113, 7057, 86813$ \\ [0.05ex]
  \hline
  \rule{0pt}{1em}\multirow{4}{*}{$4$}
  & $725$  & $16$ & $5$ & $1$ & $137\, \dots\, 9011, 13681, 14759, 35977, 264919, 794111, 1716109$\\
  & $1957$ & $4$  & $2$ &     & $47, 223, 863, 2593, 28789, 185599$ \\
  & $2777$ & $3$  & $2$ &     & $197, 359, 1301, 28057, 63607$ \\
  & $8768$ &      & $1$ &     & $10169$ \\ [0.05ex]
  \hline
  \rule{0pt}{1em}\multirow{2}{*}{$5$}
  & $70601$ & $2$ & $1$ &     & $193, 719, 26777$ \\
  & $11^4$  &     & $1$ &     & $86161$ \\ [0.05ex]
  \hline
  \rule{0pt}{1em}\multirow{1}{*}{$6$}
  & $13^5$  &   & $1$ &     & $171713$ \\ [0.05ex]
  \hline
\end{tabular}
~\\

\bibliographystyle{plain}
\bibliography{ssjbib}{}

\begin{thebibliography}{10}

\bibitem{adventures}
Sarah {Arpin}, Catalina {Camacho-Navarro}, Kristin {Lauter}, Joelle {Lim},
  Kristina {Nelson}, Travis {Scholl}, and Jana {Sot{\'a}kov{\'a}}.
\newblock {Adventures in Supersingularland}.
\newblock {\em arXiv e-prints}, page arXiv:1909.07779, September 2019.

\bibitem{bgr}
Michael~A. Bennett, Adela Gherga, and Andrew Rechnitzer.
\newblock Computing elliptic curves over {$\Bbb{Q}$}.
\newblock {\em Math. Comp.}, 88(317):1341--1390, 2019.

\bibitem{bbbccdldlrsv}
Alex~J. {Best}, Jonathan {Bober}, Andrew~R. {Booker}, Edgar {Costa}, John
  {Cremona}, Maarten {Derickx}, David {Lowry-Duda}, Min {Lee}, David {Roe},
  Andrew~V. {Sutherland }, and John {Voight}.
\newblock {Computing classical modular forms}.
\newblock {\em arXiv e-prints}, page arXiv:2002.04717, February 2020.

\bibitem{bklhpr}
Manjul Bhargava, Daniel~M. Kane, Hendrik~W. Lenstra, Jr., Bjorn Poonen, and
  Eric Rains.
\newblock Modeling the distribution of ranks, {S}elmer groups, and
  {S}hafarevich-{T}ate groups of elliptic curves.
\newblock {\em Camb. J. Math.}, 3(3):275--321, 2015.

\bibitem{bhargava-shankar}
Manjul Bhargava and Arul Shankar.
\newblock Binary quartic forms having bounded invariants, and the boundedness
  of the average rank of elliptic curves.
\newblock {\em Ann. of Math. (2)}, 181(1):191--242, 2015.

\bibitem{antwerp}
B.~J. Birch and W.~Kuyk, editors.
\newblock {\em Modular functions of one variable. {IV}}.
\newblock Lecture Notes in Mathematics, Vol. 476. Springer-Verlag, Berlin-New
  York, 1975.

\bibitem{brent-kung}
R.~P. Brent and H.~T. Kung.
\newblock Fast algorithms for manipulating formal power series.
\newblock {\em J. Assoc. Comput. Mach.}, 25(4):581--595, 1978.

\bibitem{bcdt}
Christophe Breuil, Brian Conrad, Fred Diamond, and Richard Taylor.
\newblock On the modularity of elliptic curves over {$\bold Q$}: wild 3-adic
  exercises.
\newblock {\em J. Amer. Math. Soc.}, 14(4):843--939, 2001.

\bibitem{broker}
Reinier Br\"{o}ker.
\newblock Constructing supersingular elliptic curves.
\newblock {\em J. Comb. Number Theory}, 1(3):269--273, 2009.

\bibitem{brumer-mcguinness}
Armand Brumer and Ois\'{\i}n McGuinness.
\newblock The behavior of the {M}ordell-{W}eil group of elliptic curves.
\newblock {\em Bull. Amer. Math. Soc. (N.S.)}, 23(2):375--382, 1990.

\bibitem{charles-lauter-goren}
Denis~X. Charles, Kristin~E. Lauter, and Eyal~Z. Goren.
\newblock Cryptographic hash functions from expander graphs.
\newblock {\em J. Cryptology}, 22(1):93--113, 2009.

\bibitem{cohen}
Henri Cohen.
\newblock {\em A course in computational algebraic number theory}, volume 138
  of {\em Graduate Texts in Mathematics}.
\newblock Springer-Verlag, Berlin, 1993.

\bibitem{cooley-tukey}
James~W. Cooley and John~W. Tukey.
\newblock An algorithm for the machine calculation of complex {F}ourier series.
\newblock {\em Math. Comp.}, 19:297--301, 1965.

\bibitem{clrs}
Thomas~H. Cormen, Charles~E. Leiserson, Ronald~L. Rivest, and Clifford Stein.
\newblock {\em Introduction to algorithms}.
\newblock MIT Press, Cambridge, MA, third edition, 2009.

\bibitem{cremona}
John Cremona.
\newblock ecdata: 2016-10-17, October 2016.

\bibitem{darmon}
Henri Darmon.
\newblock {\em Rational points on modular elliptic curves}, volume 101 of {\em
  CBMS Regional Conference Series in Mathematics}.
\newblock Published for the Conference Board of the Mathematical Sciences,
  Washington, DC; by the American Mathematical Society, Providence, RI, 2004.

\bibitem{dechene}
Isabelle D{\'{e}}ch{\`{e}}ne.
\newblock Quaternion algebras and the graph method for elliptic curves, 1998.
\newblock Thesis (MSc.)--McGill University (Canada),
  \url{https://www.math.mcgill.ca/darmon/theses/dechene-master/thesis.pdf}.

\bibitem{diamond-shurman}
Fred Diamond and Jerry Shurman.
\newblock {\em A first course in modular forms}, volume 228 of {\em Graduate
  Texts in Mathematics}.
\newblock Springer-Verlag, New York, 2005.

\bibitem{ehlmp}
Kirsten Eisentr\"{a}ger, Sean Hallgren, Kristin Lauter, Travis Morrison, and
  Christophe Petit.
\newblock Supersingular isogeny graphs and endomorphism rings: reductions and
  solutions.
\newblock In {\em Advances in cryptology---{EUROCRYPT} 2018. {P}art {III}},
  volume 10822 of {\em Lecture Notes in Comput. Sci.}, pages 329--368.
  Springer, Cham, 2018.

\bibitem{gross}
Benedict~H. Gross.
\newblock Heights and the special values of {$L$}-series.
\newblock In {\em Number theory ({M}ontreal, {Q}ue., 1985)}, volume~7 of {\em
  CMS Conf. Proc.}, pages 115--187. Amer. Math. Soc., Providence, RI, 1987.

\bibitem{harron-snowden}
Robert Harron and Andrew Snowden.
\newblock Counting elliptic curves with prescribed torsion.
\newblock {\em J. Reine Angew. Math.}, 729:151--170, 2017.

\bibitem{jao-defeo}
David Jao and Luca De~Feo.
\newblock Towards quantum-resistant cryptosystems from supersingular elliptic
  curve isogenies.
\newblock In {\em Post-quantum cryptography}, volume 7071 of {\em Lecture Notes
  in Comput. Sci.}, pages 19--34. Springer, Heidelberg, 2011.

\bibitem{kohel}
David~Russell Kohel.
\newblock {\em Endomorphism rings of elliptic curves over finite fields}.
\newblock ProQuest LLC, Ann Arbor, MI, 1996.
\newblock Thesis (Ph.D.)--University of California, Berkeley.

\bibitem{lmfdb}
The {LMFDB Collaboration}.
\newblock The {L}-functions and modular forms database.
\newblock \url{http://www.lmfdb.org}, 2020.
\newblock [Online; accessed 20 September 2020].

\bibitem{martin_alsize}
Kimball Martin.
\newblock Refined dimensions of cusp forms, and equidistribution and bias of
  signs.
\newblock {\em J. Number Theory}, 188:1--17, 2018.

\bibitem{martin}
Kimball Martin.
\newblock An on-average maeda-type conjecture in the level aspect.
\newblock {\em Proceedings of the American Mathematical Society}, page~1, Sep
  2020.

\bibitem{mestre}
J.-F. Mestre.
\newblock La m\'{e}thode des graphes. {E}xemples et applications.
\newblock In {\em Proceedings of the international conference on class numbers
  and fundamental units of algebraic number fields ({K}atata, 1986)}, pages
  217--242. Nagoya Univ., Nagoya, 1986.

\bibitem{mills}
W.~H. Mills.
\newblock Continued fractions and linear recurrences.
\newblock {\em Math. Comp.}, 29:173--180, 1975.

\bibitem{ppvw}
Jennifer Park, Bjorn Poonen, John Voight, and Melanie~Matchett Wood.
\newblock A heuristic for boundedness of ranks of elliptic curves.
\newblock {\em J. Eur. Math. Soc. (JEMS)}, 21(9):2859--2903, 2019.

\bibitem{poonen}
Bjorn Poonen.
\newblock Heuristics for the arithmetic of elliptic curves.
\newblock In {\em Proceedings of the {I}nternational {C}ongress of
  {M}athematicians---{R}io de {J}aneiro 2018. {V}ol. {II}. {I}nvited lectures},
  pages 399--414. World Sci. Publ., Hackensack, NJ, 2018.

\bibitem{rabin}
Michael~O. Rabin.
\newblock Probabilistic algorithms in finite fields.
\newblock {\em SIAM J. Comput.}, 9(2):273--280, 1980.

\bibitem{silverman}
Joseph~H. Silverman.
\newblock {\em The arithmetic of elliptic curves}, volume 106 of {\em Graduate
  Texts in Mathematics}.
\newblock Springer, Dordrecht, second edition, 2009.

\bibitem{mse}
Math Stackexchange.
\newblock Formulas for the (top) coefficients of the characteristic polynomial
  of a matrix.
\newblock \url{https://math.stackexchange.com/questions/23899/}.
\newblock [Online; accessed 20 September 2020].

\bibitem{sturm}
Jacob Sturm.
\newblock On the congruence of modular forms.
\newblock In {\em Number theory ({N}ew {Y}ork, 1984--1985)}, volume 1240 of
  {\em Lecture Notes in Math.}, pages 275--280. Springer, Berlin, 1987.

\bibitem{sutherland}
Andrew~V. Sutherland.
\newblock Isogeny volcanoes.
\newblock In {\em A{NTS} {X}---{P}roceedings of the {T}enth {A}lgorithmic
  {N}umber {T}heory {S}ymposium}, volume~1 of {\em Open Book Ser.}, pages
  507--530. Math. Sci. Publ., Berkeley, CA, 2013.

\bibitem{taylor-wiles}
Richard Taylor and Andrew Wiles.
\newblock Ring-theoretic properties of certain {H}ecke algebras.
\newblock {\em Ann. of Math. (2)}, 141(3):553--572, 1995.

\bibitem{watkins}
Mark Watkins.
\newblock Some heuristics about elliptic curves.
\newblock {\em Experiment. Math.}, 17(1):105--125, 2008.

\bibitem{wiedemann}
Douglas~H. Wiedemann.
\newblock Solving sparse linear equations over finite fields.
\newblock {\em IEEE Trans. Inform. Theory}, 32(1):54--62, 1986.

\bibitem{wiles}
Andrew Wiles.
\newblock Modular elliptic curves and {F}ermat's last theorem.
\newblock {\em Ann. of Math. (2)}, 141(3):443--551, 1995.

\end{thebibliography}

\end{document}